\documentclass[reqno,11pt]{amsart}
\usepackage{amsthm,amsfonts,amssymb,euscript,mathrsfs,graphics,color,amsmath,latexsym,marginnote,hyperref}
\usepackage[latin1]{inputenc} 
\usepackage{amsaddr}

\theoremstyle{plain}

\usepackage{times}
\numberwithin{equation}{section}

\newtheorem{Teo}{Theorem}
\newtheorem{Def}{Definition}
\newtheorem{theorem}{Theorem}[section]

\newtheorem{remark}[theorem]{Remark}
\newtheorem{remarks}[theorem]{Remark}
\newtheorem{definition}[theorem]{Definition}
\newcommand{\be}{\begin{equation}}
\newcommand{\ee}{\end{equation}}

\newcommand{\e}{\varepsilon}

\newcommand{\R}{\mathbb R}

\newcommand{\Z}{\mathbb Z}

\newcommand{\N}{\mathbb N}
\newcommand{\T}{\mathbb T}
\newcommand{\sign}{{\rm sign}}

\newcommand{\s }{\sigma }
\newcommand{\ii }{{\rm i} }

\newcommand{\vphi}{\varphi }

\newcommand{\dps}{\displaystyle}

\newcommand{\pa}{\partial}

\def\bar{\overline}

\def\cal{\mathcal}

\def\ba{\begin{aligned}}
\def\ea{\end{aligned}}
\def\beginm{\begin{multline}}
\def\endm{\end{multline}}

 \providecommand{\sign}{\mathrm{sgn}\,}

\makeatletter

\def\l@subsection{\@tocline{2}{0pt}{2.5pc}{5pc}{}}
\def\l@subsubsection{\@tocline{3}{0pt}{4.5pc}{5pc}{}}

\renewcommand\tocchapter[3]{%
  \indentlabel{\@ifnotempty{#2}{\ignorespaces#2.\quad}}#3%
}
\newcommand\@dotsep{4.5}
\def\@tocline#1#2#3#4#5#6#7{\relax
  \ifnum #1>\c@tocdepth 
  \else
    \par \addpenalty\@secpenalty\addvspace{#2}%
    \begingroup \hyphenpenalty\@M
    \@ifempty{#4}{%
      \@tempdima\csname r@tocindent\number#1\endcsname\relax
    }{%
      \@tempdima#4\relax
    }%
    \parindent\z@ \leftskip#3\relax \advance\leftskip\@tempdima\relax
    \rightskip\@pnumwidth plus1em \parfillskip-\@pnumwidth
    #5\leavevmode\hskip-\@tempdima{#6}\nobreak
    \leaders\hbox{$\m@th\mkern \@dotsep mu\hbox{.}\mkern \@dotsep mu$}\hfill
    \nobreak
    \hbox to\@pnumwidth{\@tocpagenum{#7}}\par
    \nobreak
    \endgroup
  \fi}
\def\l@subsection{\@tocline{2}{0pt}{2.5pc}{5pc}{}}

\begin{document}

\title{Time quasi-periodic traveling gravity water waves in infinite depth}

\date{}

\author{Roberto Feola}
\address{University of Nantes,\\
e-mail: roberto.feola@univ-nantes.fr}

\author{Filippo Giuliani}
\address{UPC, Barcelona, \\e-mail: filippo.giuliani@upc.edu}


 \begin{abstract} 
We present the recent result \cite{FGww} concerning the existence of quasi-periodic in time traveling waves for the $2d$ 
pure gravity water waves system in infinite depth. We provide the first existence result of quasi-periodic water waves solutions bifurcating from a completely resonant elliptic fixed point. The proof is based on a Nash-Moser scheme, Birkhoff normal form methods and pseudo-differential calculus techniques. We deal with the combined problems 
of small divisors
and the fully-nonlinear nature of the equations.
\end{abstract}

\maketitle

\section{Introduction}

The aim of this note is to present the results obtained in \cite{FGww}
in which we prove the existence and the linear stability of small amplitude, quasi-periodic traveling solutions
for the $2d$ 
pure gravity water waves system in infinite depth. 
By a quasi-periodic traveling wave we mean a motion that, 
at the first orders of amplitude, is a superposition 
of an arbitrarily large number of periodic traveling waves with rationally independent frequencies.
The irrationality of the frequencies of oscillations 
excludes the existence of a moving frame for 
which such motions are stationary. The existence of quasi-periodic in time water waves have been proved only recently, we refer to the papers \cite{BM1}, \cite{BBHM}. In these works the authors deal with cases in which the linear frequencies of oscillations are modulated by physical parameters, like the capillarity of the fluid or the depth of the ocean. 
In the pure gravity case with infinite depth the lack of such parameters
allows the existence of infinite-dimensional subspaces which are invariant for the linearized problem and are filled by \emph{periodic} in time solutions. We refer to this situation as a \emph{completely resonant} case (see also \cite{IoossPloTol1}). Usually the complete resonance is referred as a stronger property of the linearized equation: all the solutions are periodic. Nevertheless the difficulties that we encounter are the same, indeed our solutions bifurcate from (possibly) periodic solutions. Then a careful nonlinear analysis is required to extract, directly from the system, some parameters that tune in an efficient way the frequencies of the expected solutions.
This task is complicated by the presence of \emph{small divisors}
and the \emph{fully-nonlinear} nature of the equations.

\subsection{Formulation of the problem}
We consider an incompressible and irrotational perfect fluid, under the action of gravity
occupying, at time
$t$, a two dimensional domain with infinite depth, periodic in the horizontal variable, given by 
\begin{equation*}
    {\mathcal D}_{\eta} := \big\{ (x,y)\in \T \times\R \, ;  \ - \infty <y<\eta(t,x) \big\},  \quad 
    \T := \R \slash (2 \pi \Z) \, .
\end{equation*}
The velocity field in the time dependent domain 
$ {\mathcal D}_{\eta} $ is the gradient of a harmonic function $\Phi$, called the velocity potential.
The time-evolution of the fluid is determined by a system of equations for 
the free surface  $\eta(t,x) $, 
and the function $\psi(t,x) = \Phi(t,x,\eta(t,x))$ which is
 the restriction 
of the velocity potential $\Phi$ to the free interface.
Given the shape $\eta(t,x)$ of the domain $ {\cal D}_\eta $ 
and the Dirichlet value $\psi(t,x)$ of the velocity potential at the top boundary, one can recover 
$\Phi(t,x,y)$ as the unique solution of the elliptic problem 
\begin{equation} \label{BoundaryPr}
\Delta \Phi = 0  \ \text{in } 
{\cal D}_\eta  \, , \quad 
\partial_y \Phi \to 0  \ \text{as } y \to - \infty \, , \quad 
\Phi = \psi \ \;\; \text{on }\;\; \{y = \eta(t,x)\}. 
\end{equation}

According to Zakharov \cite{Zak1} and Craig-Sulem \cite{CrSu}
the $(\eta,\psi)$ variables satisfy the gravity water waves system
\begin{equation} \label{eq:113}
\begin{cases}
    \partial_t \eta = G(\eta)\psi \cr
\partial_t\psi = \dps -g\eta  -\frac{1}{2} \psi_x^2 +  \frac{1}{2}\frac{(\eta_x  \psi_x + G(\eta)\psi)^2}{1+\eta_x^2}
\end{cases}
\end{equation}
where $ G(\eta)\psi $ is the Dirichlet-Neumann operator  
\begin{equation*}
  G(\eta)\psi :=  (\partial_y\Phi -\eta_x \partial_x\Phi)(t,x,\eta(t,x))\,.
\end{equation*}
Without loss of generality, we set the gravity constant to $g = 1$.

It was first observed by Zakharov \cite{Zak1} that \eqref{eq:113} is a Hamiltonian system with respect to the symplectic form $d\psi\wedge d\eta$ and it can be written as
\begin{equation}\label{HS}
\begin{aligned}
& \qquad \pa_t \eta = \nabla_\psi H (\eta, \psi) \, , \quad  \pa_t \psi = - \nabla_\eta H (\eta, \psi)  \, , 
\end{aligned}
\end{equation}
where $ \nabla $ denotes the $ L^2 $-gradient, with Hamiltonian
\begin{equation}\label{Hamiltonian}
H(\eta, \psi) := \frac12 \int_\T \psi \, G(\eta ) \psi \, dx + \frac{1}{2} \int_{\T} \eta^2  \, dx
\end{equation}
given by the sum of the kinetic 
and potential energy of the fluid.
The invariance of the system \eqref{eq:113} in the $y$ and $x$ variable implies the existence of two prime integrals, respectively the ``mass'' $\int_\T \eta \, dx$ and the momentum 
\begin{equation}\label{Hammomento}
\mathtt{M}:= \int_{\T} \eta_x (x) \psi (x)  \, dx. 
\end{equation}
Since  the ocean has infinite depth,   if $\Phi$ solves \eqref{BoundaryPr},
then $\Phi_{c}(x,y):=\Phi(x,y-c)$ solves the same problem in $\mathcal{D}_{\eta+c}$ 
assuming the Dirichlet datum $\psi$ at the free boundary
 $\eta+c $. Therefore $G(\eta+c)=G(\eta) $, for all $ c \in \R $, 
and $ \int_\T \nabla_\eta K \, dx = 0$ where $ K := \frac12 \int_{\T} \psi\, G(\eta) \psi \, dx $ denotes the kinetic energy.
Then
$ \widehat \eta_0 (t)  := \frac{1}{2\pi} \int_{\T} \eta (t, x) \, dx $, 
$\widehat \psi_0 (t) := \frac{1}{2\pi} \int_{\T} \psi (t, x) \, dx $
evolve according to the decoupled equations
\be\label{medie00}
\pa_t {\widehat \eta}_0 (t) = 0 \, , \quad \pa_t  \widehat \psi_0 (t) = - g \widehat \eta_0 (t)  .
\ee
Hence
we may restrict the study of the dynamics to the invariant subspace of $(\eta, \psi)$ such that
\begin{equation}\label{zeromode}
\int_{\T} \eta \, dx = \int_{\T} \psi \, dx = 0  \, . 
\end{equation}

The Hamiltonian \eqref{Hamiltonian}  is defined on the spaces
\begin{equation*}
(\eta, \psi) \in H^s_0 (\T; \mathbb{R}) \times {H}_0^s (\T; \mathbb{R})  
\end{equation*}
where $ H^s (\T;\mathbb{R})$, $ s\in \R $, denotes the Sobolev space of $ 2 \pi $-periodic functions of $ x $, 
and $H^s_0(\T;\mathbb{R})$ is the subspace of 
$H^s(\T;\mathbb{R})$ of zero average functions. 

\smallskip

Small amplitude solutions are close to the solutions of the linearized system of \eqref{eq:113}
at the equilibrium $(\eta,\psi)=(0,0)$, namely
\begin{equation} \label{eq:113linear}
    \partial_t \eta =G(0)\psi\,,\qquad 
\partial_t\psi = \dps -\eta  
\end{equation}
where 
 the Dirichlet-Neumann operator 
at the flat surface $\eta=0$ is the Fourier multiplier
$G(0)=|D|$.
Passing to the complex coordinates, the system \eqref{eq:113linear}
is equivalent to the following equation
\begin{equation} \label{eq:113CWW}
\pa_{t}u=-\ii |D|^{\frac{1}{2}} u\,,\qquad 
u=\frac{1}{\sqrt{2}}\big( |D|^{-\frac{1}{4}}\eta+\ii|D|^{\frac{1}{4}}\psi \big)\,.
\end{equation}
 The linear solutions have the form
\begin{equation*}
\begin{aligned}
u(t,x)&=\sum_{j\in \mathbb{Z}\setminus\{0\}}
u_{j}(0)e^{-\ii\sqrt{|j|}t+\ii jx}\,,\\
 u_j(0)&:=\frac{1}{\sqrt{2}}\big( |j|^{-\frac{1}{4}}\eta_j(0)+\ii|j|^{\frac{1}{4}}\psi_{j}(0) \big)\,.
\end{aligned}
\end{equation*}
Such solutions  can
 be either periodic, quasi-periodic or almost-periodic depending on the Fourier support.
 We remark that solutions initially Fourier supported on the infinite-dimensional invariant subspace $\{ e^{\mathrm{i} n^2 x} : n\in\mathbb{Z} \}$ are periodic in time, hence we say that the equation \eqref{eq:113CWW} is \emph{completely resonant}.
The frequency of oscillation of the $j$-th mode is $\sqrt{|j|}$ and
we refer to the map 
\begin{equation}\label{dispersionLaw}
j\to \sqrt{|j|}\,,\quad j\in \mathbb{Z}\setminus\{0\}\,,
\end{equation}
as the \emph{dispersion law} of \eqref{eq:113linear}.
We note that the dispersion law is even in $j$, hence there are infinitely many multiple eigenvalues, and it grows \emph{sublinearly}.

\subsection{Main results}
We define quasi-periodic traveling waves in the following way.
\begin{Def}{\bf (Quasi-periodic traveling waves).}
$(i)$ We say that a function 
\[
(\eta(t,x),\psi(t,x)) : \mathbb{R}\times \mathbb{T}\to\mathbb{R}^{2}
\]
is a \emph{quasi}-periodic solution of \eqref{eq:113} with  irrational 
frequency vector $\omega\in \mathbb{R}^{\nu}$, if there is an embedding 
\begin{equation}\label{6.6}
\begin{aligned}
\T^{\nu} & \to  
H_0^1(\mathbb{T};\mathbb{R})\times
 {H}_0^1(\mathbb{T};\mathbb{R})\\
 \varphi & \mapsto U(\varphi,x):=(\tilde{\eta}(\varphi, x), \tilde{\psi}(\varphi, x))
 \end{aligned}
\end{equation}
such that $(\eta(t,x),\psi(t,x))=U(\omega t,x)$ solves \eqref{eq:113}.

\noindent
$(ii)$ A quasi-periodic solution is \emph{traveling} with \emph{velocity vector}
$\mathtt{v}\in \mathbb{Z}^{\nu}$ if there is a function 
$\widetilde{U} : \mathbb{T}^{\nu}\to \mathbb{R}^{2}$ such that
\begin{equation}\label{travellingWaves}
(\eta(t,x),\psi(t,x))=U(\omega t, x)=\widetilde{U}(\omega t-\mathtt{v}x)\,.
\end{equation}
\end{Def}
\begin{remark}
We remark that an 
embedding $U$ satisfies  \eqref{travellingWaves} if and only if $U(\vphi,x)$
solves the transport equation 
\begin{equation*}
\mathtt{v}\cdot \partial_{\varphi} U+X_{\mathtt{M}}(U)=0\,, 
\qquad X_{\mathtt{M}}(U)=(\tilde{\eta}_x, \tilde{\psi}_x)\,,
\end{equation*}
where $\mathtt{M}$ is the momentum Hamiltonian in \eqref{Hammomento}.
\end{remark}
We shall construct 
such solutions localized  
 in Fourier space at 
$\nu$ distinct  \textit{tangential sites} $S:=S^{+}\cup S^{-}$
\begin{equation*}
S^+:=\{\overline{\jmath}_1, \dots, \overline{\jmath}_m \}\subset \mathbb{N}\setminus\{0\}\,,
\quad
S^-:=\{ \overline{\jmath}_{m+1}, \dots, \overline{\jmath}_{\nu}\}\subset -\mathbb{N}\setminus\{0\}\,,
\end{equation*}
for some  $1\leq m\leq \nu$ and where
\begin{equation*}
k \neq -j\,,\quad \forall j\in S^+,\,\,\forall k\in S^-.
\end{equation*}
The solutions of \eqref{eq:113linear} that originate by exciting the tangential modes are superpositions of periodic traveling linear waves with velocity $\overline{\jmath}_i$ and frequency $\sqrt{|\overline{\jmath}_i|}$. Such motions are quasi-periodic (or periodic) traveling waves of the form \eqref{travellingWaves} with frequency vector
\begin{equation*}
\overline{\omega}:=\left(
\sqrt{|\overline{\jmath}_1|},\ldots, \sqrt{|\overline{\jmath}_{\nu}|}
\right)\in\mathbb{R}^{\nu}
\end{equation*}
and velocity vector
\begin{equation}\label{velocityvec}
\mathtt{v}:=\left( \overline{\jmath}_1, \dots, \overline{\jmath}_{\nu} \right)\in\mathbb{Z}^{\nu}.
\end{equation}
We construct small amplitude, quasi-periodic traveling waves solutions of \eqref{eq:113} 
which are ``close'' to the linear ones, namely
they will be of the form
\begin{equation}\label{SoluzioneEsplicitaWW}
\begin{aligned}
\eta(t,x)&=\sum_{j\in S}
\sqrt{2\zeta_{j}}|j|^{\frac{1}{4}}\Big(\cos(\omega_{j}t)\cos(j x)+
\sin(\omega_{j}t)\sin(j x)\Big)
+o(\sqrt{|\zeta|})\,,\\
\psi(t,x)&=\sum_{j\in S}
\sqrt{2\zeta_{j}}
|j|^{-\frac{1}{4}}\Big(
\cos(\omega_{j}t)\sin(jx)-
\sin(\omega_{j}t)\cos(jx)\Big)+o(\sqrt{|\zeta|})\,, 
\end{aligned}
\end{equation}
%
where $\zeta=(\zeta_j)_{j\in S}$ with $\zeta_j>0$, $\omega=\bar\omega +O(\lvert \zeta \rvert)$ and $o(\sqrt{\lvert \zeta \rvert})$
is meant  in the 
 $H^{s}$-topology with $s$ large. The vectors $(\sqrt{\zeta_j})_{j\in S}$ represent the \emph{amplitudes} of the approximate solution from which we have the bifurcation.
 Our main result holds for a suitable choice of the tangental sites 
 $\overline{\jmath}_i$, that  we prove to be   \emph{generic}. When we refer to a generic choice of the tangential sites we mean that the $\overline{\jmath}_i$'s are chosen such that the vector $\left( \overline{\jmath}_1, \dots, \overline{\jmath}_{\nu} \right)$ is not a zero of a certain non-trivial polynomial $\mathbb{C}^{\nu}\to\mathbb{C}$. We remark that such a choice is equivalent to the choice of the velocity vector $\mathtt{v}$ in \eqref{velocityvec}.

Denoting by $B(0,\varrho)$ 
the ball centered at the origin of $\mathbb{R}^{\nu}$ 
of radius $\varrho>0$, 
our result can be stated as follows.
\begin{Teo}{\bf (Quasi-periodic traveling  gravity waves).}\label{thm:main}
Let $\nu\geq 1$.
For a generic
choice of the velocity vector $\mathtt{v}$ as in \eqref{velocityvec} there exist $s\gg1$,  $0<\varrho\ll1$     
and a positive measure  
 Cantor-like set $\mathfrak{A}\subseteq B(0,\varrho)$
 such that
the following holds. 
For any $\zeta\in \mathfrak{A}$,
the equation \eqref{eq:113} possesses a small amplitude, linearly stable, 
quasi-periodic solution $(\eta, \psi)(t,x; \zeta)=U(\omega t,x; \zeta)$ of the form \eqref{SoluzioneEsplicitaWW} which is a traveling wave with velocity vector $\mathtt{v}$, 
$U(\vphi,x)\in H^{s}(\mathbb{T}^{\nu+1},\mathbb{R}^{2})$ and 
$\omega:=\omega(\zeta)\in\mathbb{R}^{\nu}$ is a diophantine frequency vector. Moreover for  $0<\e\le \sqrt{\varrho}$, the set  $\mathfrak A$ has asymptotically full relative measure in $[\e^2,2\e^2]^\nu$.
\end{Teo}

To the best of our knowledge the above theorem is the first existence result concerning
quasi-periodic solutions of the water waves equations bifurcating from a completely resonant elliptic fixed point. 

\smallskip

\noindent
We mention Iooss-Plotnikov-Toland  \cite{IoossPloTol1} and Iooss-Plotnikov \cite{IossPlo111} which proved for the completely resonant equation  \eqref{eq:113} the existence of \emph{standing periodic} solutions.\\
Concerning the traveling waves, we cite Craig-Nicholls \cite{CN1} which proved the existence of \emph{periodic} solutions 
in
 the gravity-capillary case with space periodic boundary conditions.
In absence of capillarity the existence of periodic traveling waves is 
a small divisor problem. 
This case has been treated by 
 Iooss and Plotnikov in \cite{IoossPlo1, IoossPlo2}.\\
 The first results on the existence of quasi-periodic waves are quite recent  
and are due to Berti-Montalto \cite{BM1} for the gravity-capillary 
case with infinite depth and by Baldi-Berti-Haus-Montalto \cite{BBHM} 
for the pure gravity case with finite depth for $2$d oceans. 
In both cases the existence of quasi-periodic solutions
is provided for some asymptotically full-measure set of the parameters of the problem, respectively capillarity and depth (or equivalently wavelength). \\
As far as we know, all the previous results on periodic and quasi-periodic in time water waves take advantage from the presence of physical parameters and / or assumptions of parity conditions.\\
The purpose of the paper \cite{FGww} is to address two natural questions:\\
 (i) we work on a fixed equation for which the only possible parameters to modulate are the initial data of the solutions.\\
(ii) we look for a general class of quasi-periodic traveling waves which are free from restrictions of parity in the spatial and time variables.

\vspace{0.9em}
\noindent
{\bf Acknowledgements.}
The authors would like to thank the anonymous referee for very valuable and useful comments.\\
Roberto Feola has been supported of the Centre Henri Lebesgue 
ANR-11-LABX-0020-01 and by ANR-15-CE40-0001-02 
``BEKAM'' of the Agence Nationale de la Recherche.
Filippo Giuliani has received funding from 
the European Research Council (ERC) 
under the European Union's Horizon 2020
research and innovation programme under grant agreement 
No 757802.

\section{Comments of the main result}
Now we discuss the main issues and the novelties of the paper \cite{FGww}.

\smallskip
\noindent
$\bullet$  The general form of the linear frequencies of oscillations for the water waves equations is the following 
\[
\sqrt{| j | \tanh(\mathtt{h} |j|) (g +\kappa j^{2})},
\]
where $\mathtt{h}$ and $\kappa$ are respectively the depth and 
the capillarity of the fluid. 
If $\mathtt{h}<\infty$ or $\kappa\neq 0$ such parameters may be used to impose non-resonance conditions that are needed for the search of quasi-periodic solutions (see for instance \cite{BBHM}, \cite{BM1}, \cite{PloTol1}). In our case $\mathtt{h}=\infty$ and $\kappa= 0$, thus the linear frequencies of oscillations are
$
\sqrt{g \,|j|}
$
and the elements of the infinite dimensional space $\mbox{span}\{ e^{\mathrm{i} |n| t}\,e^{\mathrm{i} n^2 x} : n\in\mathbb{Z} \}$ are \emph{periodic} solutions of the linearized problem at the origin (completely resonant case).\\
The physical parameter $g$ clearly does not modulate the frequencies, hence if we look for quasi-periodic solutions
 we need to extract parameters directly from the nonlinearities of the equation. We do that by means of \emph{Birkhoff normal form} (BNF) techniques.
 In this way the bifurcation parameters are essentially the ``initial data'' 
 or the amplitudes
of an appropriate approximate solution 
 from which the bifurcation occurs (see \eqref{SoluzioneEsplicitaWW}). 
The choice of the Fourier support of such approximate solution 
plays a fundamental role in proving some non-degeneracy conditions.
Roughly speaking, both the amplitudes $\sqrt{\zeta_{\overline{\jmath}_i}}$ and the tangential sites $\overline{\jmath}_i$
will be ``parameters'' of our problem.


\smallskip
\noindent
$\bullet$
In performing BNF procedures we shall deal with resonances among linear frequencies.
It is known that the pure gravity case in infinite depth 
has no $3$-waves resonant interactions. 
On the other hand, there are many \emph{non-trivial} $4$-wave
interactions, called 
\emph{Benjamin-Feir} resonances (see \eqref{straBFR}).
 We then exploit
a fundamental property of the pure gravity waver waves Hamiltonian 
\eqref{Hamiltonian} in infinite depth: 
the formal integrability, up to order four, of the 
Birkhoff normal form.
This has been proved in \cite{Zak2}, \cite{CW}, 
by showing explicit key algebraic cancellations 
occurring for the coefficients of the Hamiltonian.

\smallskip
\noindent
$\bullet$ In order to show that the ``initial data'' of the expected 
solutions tune in an efficient way the frequencies we 
shall find the explicit expression of the first order 
corrections of the tangential frequencies and of the spectrum 
of the linearized operator in the normal directions. 
 We use an identification argument of normal 
 forms to detect them, in the spirit of \cite{FGP1}, based on the presence of approximate constants 
 of motion.
Actually, even after this procedure, the 
Hamiltonian is still partially degenerate. Indeed it turns out that there are 
 a finite number of eigenvalues which are still in resonance.
We overcome this difficulty 
by passing to suitable rotating coordinates.

\smallskip
\noindent
$\bullet$ 
We exhibit the existence of a wide class of traveling quasi-periodic solutions with no parity restrictions in time and space by using the Hamiltonian structure and the $x$-translation invariance of \eqref{eq:113}. It is well known that 
the water waves system \eqref{eq:113} 
exhibits additional  symmetries.
For instance the vector field $X_{H}$
in \eqref{HS} is \begin{itemize}
\item[(i)]  \emph{reversible} with respect to  the involution
 \begin{equation*}
S : \left(\begin{matrix} 
\eta(x) \\ \psi(x)
\end{matrix}\right)\;  \mapsto \;
 \left(\begin{matrix} 
\eta(-x) \\ -\psi(-x)
\end{matrix}\right)
\end{equation*}
i.e. it satisfies $X_{H}\circ S=-S\circ X_{H}$;
\item[(ii)]
 \emph{even-to-even}, 
i.e. maps $x$-even functions into $x$-even functions. 
\end{itemize}

In several papers (for instance \cite{BBHM}, \cite{BM1}, \cite{IoossPloTol1}, \cite{PloTol1}) such symmetries are adopted to remove degeneracies due to translation invariance in $x$ and $t$.
Since we do \emph{not} look for solutions in the subspace of reversible functions, namely
\begin{equation*}
\{(\eta, \psi) \quad \mbox{such that} \quad(\eta, \psi)(t, x)
=(\eta, -\psi)(-t, -x)=S(\eta, \psi)(-t, x)\},
\end{equation*}
the existence of a solution $U(t, x)$ of \eqref{eq:113} implies the existence of a possibly different solution $S U(-t, x)$.

\smallskip
\noindent
$\bullet$ The vector field in \eqref{eq:113} is a \emph{singular perturbation}
of the linearized  system at the origin \eqref{eq:113linear}
since
the nonlinearity
contains derivatives of the first order, while \eqref{eq:113CWW}
contains only derivatives of order $1/2$.
We remark that this is not the case 
when the capillarity $\kappa\neq 0$.

\smallskip
\noindent
$\bullet$ The dispersion law \eqref{dispersionLaw} is \emph{sub-linear}.
This is a major difference in developing KAM theory 
for equations with \emph{super-linear} dispersion law, such as in the gravity-capillary case. 
Indeed
 a weaker dispersion law
implies \emph{bad} separation properties of the eigenvalues. The main issue concerns the verification of non-resonance conditions between the tangential frequencies and the difference of the normal ones, called ``second order Melnikov conditions''. 
%

\section{Ideas of the proof}

  Here we 
 summarize the steps of the proof of Theorem \ref{thm:main} highlighting the key ingredients.

\smallskip
\noindent
 \emph{Nash-Moser theorem of hypothetical conjugation.}
We formulate the problem of finding quasi-periodic solutions as the problem of the search for zeros of a nonlinear functional equation
\[
\mathcal{F}(\omega, U(\omega; \varphi, x))=0
\]
where $U(\omega; \varphi, x)$ is a smooth embedding of a $\nu$-dimensional torus supporting a quasi-periodic motion with frequency vector $\omega\in\mathbb{R}^{\nu}$ (see \eqref{6.6}). We implement a Nash-Moser scheme which provides the zeros of such a functional 
as limit of a sequence $(U_n)_{n\geq 0}$ of approximate solutions
\[
U_{n+1}:=U_n-\Pi_n \,\big(d \mathcal{F}(\omega; U_n)\big)^{-1}\,\Pi_n \,\mathcal{F}(U_n),
\] 
where $\Pi_n$ is a smoothing operator,
that converges in some very regular Sobolev space.
The main issues concern the invertibility of the 
linearized operator (with suitable estimates) in a neighborhood of the 
equilibrium and the search for a good 
approximate solution $U_0$ that initializes the scheme. This task is doable only for appropriate choices of the parameter frequency $\omega$. Then we prove that the set of such \emph{good} frequencies has positive Lebesgue measure.

\smallskip
\noindent
 \emph{Bifurcation analysis.}
We find the first nonlinear approximation $U_0$ as a solution of a "simplified" Hamiltonian system whose dynamics is close to \eqref{Hamiltonian}, at least for a certain time range.
Since we work in a neighborhood of an elliptic fixed point, this Hamiltonian can be obtained through a
 \emph{weak} version of the Birkhoff normal form method (WBNF), where the adjective weak refer to the fact that we just "partially" normalize the Hamiltonian \eqref{Hamiltonian}. 
More precisely we construct
a map, which is close to the identity up to a 
finite rank operator, such that 
the Hamiltonian $H$ in \eqref{Hamiltonian}, in the complex coordinates \eqref{eq:113CWW},
assumes the form $H_{\rm Birk}+R$, 
where $R$ can be considered as a small remainder in a neighborhood of the origin,
\begin{align}
H_{\rm Birk}&=H^{(2)}_{\mathbb{C}}+\mathcal{H}_{\mathbb{C}}^{(4,0)}+\mathcal{H}_{\mathbb{C}}^{(\geq 6, 0)}, \quad H^{(2)}_{\mathbb{C}}
=\sum_{j\in\mathbb{Z}} \sqrt{|j|} |u_k|^2\label{theoBirH}\\
\mathcal{H}_{\mathbb{C}}^{(4,0)} &:=
\frac{1}{4 \pi} \sum_{k \in S} |k|^3  
  |u_k|^4 
  +  \frac{1}{\pi} \sum_{\substack{k_1, k_2 \in S, \, 
  \sign(k_1) = \sign( k_2 )  \\ |k_2| < |k_1|}} |k_1| |k_2|^2  
 |u_{k_1}|^2  |u_{k_2}|^2 \,\nonumber
\end{align}
and $\mathcal{H}_{\mathbb{C}}^{(\geq 6, 0)}$ is some homogeneous integrable Hamiltonian, in the sense that it depends only on the actions $|u_k|^2$, of degree $6$. We remark that $\mathcal{H}_{\mathbb{C}}^{(4, 0)}, \mathcal{H}_{\mathbb{C}}^{(\geq 6, 0)}$ are supported on the tangential set $S$, while $R$ contains monomials supported at least on one mode out of $S$. The finite dimensional subspace
$U_S:=\{ u_k=0\,\,\,\forall k\notin S  \}$ turns out to be invariant for $H_{\rm Birk}$. We introduce the following action-angle variables on $U_S$
\[
u_k=\sqrt{I_k}\,e^{-\mathrm{i} \theta_k} \quad k\in S, \qquad u_k=z_k \quad k\notin S. 
\]
In these coordinates $\mathcal{H}_{\mathbb{C}}^{(4,0)}=\frac{1}{2} \mathbb{A} I\cdot I$, where $\mathbb{A}$ is a symmetric $\nu\times \nu$ matrix.
We prove that the Hamiltonian restricted to $U_S$ is integrable and non-degenerate
 in the sense that
the \emph{frequency-to-amplitude} map 
\[
\omega(I)=\bar{\omega}+\mathbb{A} I+O(I^2)
\] 
is a local diffeomorphism. Then we can select $U_0$ among the tori $\{I=\zeta \}$, where $\zeta$ is some fixed vector with positive components. In order to work in a neighborhood of these tori is convenient to rescale the unperturbed actions $\zeta \to \varepsilon^2 \zeta$ with a small parameter $\varepsilon>0$. The non-degeneracy of the frequency-amplitude map allows to impose non-resonance conditions such as
\begin{equation}\label{gammaDP}
|\omega\cdot \ell|\geq \frac{\gamma}{\langle \ell \rangle^{\tau}} \qquad \forall \ell\in\mathbb{Z}^{\nu} \quad \tau>\nu-1
\end{equation}
where 
\begin{equation}\label{eps2}
\gamma=o(\varepsilon^{2}).
\end{equation}
Notice that the condition \eqref{gammaDP} could be obtained just by choosing
generically the  tangential sites in $S$. However, 
 in order to 
impose higher order non-resonance conditions  required by the Nash-Moser scheme, we truly need the non-degeneracy of the frequency-to-amplitude map.
We remark that the dependence of $\gamma$ in \eqref{eps2} 
respect to $\varepsilon$ 
is a peculiarity of resonant cases, indeed this is due to the closeness of 
$\omega$ to the resonant vector $\bar{\omega}$. For non-resonant cases $\gamma$ can be considered as a small number independent of $\varepsilon$. In the search for small amplitude solutions for resonant equations this fact produces several difficulties in the bifurcation analysis. For instance, in our case, we have to implement several steps in the Birkhoff normal form procedures that we perform along the Nash-Moser iteration.\\
In the normal form analysis we have to deal 
with waves resonant interactions such as
\begin{equation*}\label{rambo3}
\sum_{i=1}^{n}\s_{i} j_i=0\,,\qquad 
\sum_{i=1}^{n}\s_{i}\sqrt{|j_{i}|}=0\,,\qquad \s_{i}=\pm\,,\;\;\;
i=1,\ldots,n\,.
\end{equation*}
We say that a $n$-tuple $(j_1,\ldots,j_{n})$  is a \emph{trivial} 
resonance
if $n$ is even, 
$\s_{i}=-\s_{i+1}$, $i=1,\ldots,n-1$ (up to permutations), and the $n$-tuple has the form
$(j,j,k,k,\ldots)$. 
It is easy to note that monomials $u_{j_1}\dots u_{j_n}$ supported on 
trivial resonances are \emph{integrable}, meaning that depend only on the actions $|u_j|^2$.
For $n=4$ there are infinitely many non-trivial resonances, called Benjamin-Feir resonances, 
which consist in 
the two parameter family of solutions
\begin{equation}\label{straBFR}
\bigcup_{\lambda\in\Z \setminus \{0\}, b\in\N} 
\Big\{ 
j_1 = -\lambda b^2, 
 \, j_2 = \lambda (b+1)^2 \, , \, j_3 = \lambda (b^2+b+1)^2, \, j_4= \lambda (b+1)^2 b^2 \Big\} \, ,
\end{equation}
with $\sigma_1=\s_3=-\sigma_2=-\s_4$. 
In  \cite{Zak2}, \cite{CW} it has been proved that 
 the coefficients of the normalized Hamiltonian 
 (obtained by a full Birkhoff normal form procedure) at order four
of the monomials corresponding to the Benjamin-Feir resonances vanish.
By using suitable algebraic arguments, we actually prove that 
such cancellations of \cite{Zak2}, \cite{CW}
occur also preforming the \emph{weak} version of BNF which involves only \emph{finitely}
many tangential sites. 
In this way we conclude the integrability
of the weak BNF Hamiltonian at degree 4 and we obtain its explicit formula 
\eqref{theoBirH}. 
In order to deal with higher order resonances we use a genericity
argument.

\smallskip
\noindent
 \emph{Normal form identification.} The frequency of the expected quasi-periodic solutions are small corrections of the linear frequencies of oscillations $\overline{\omega}$. Then the first order corrections are fundamental to impose the non-resonance conditions. Usually such corrections are obtained through a Birkhoff normal form procedure that normalizes the \emph{full} Hamiltonian at order four, i.e. that normalizes all the monomials of order four. This procedure has been implemented in \cite{CW}, \cite{Zak2} \emph{at formal level}, where the integrability of the normal form is provided by explicit computations of the coefficients of the Hamiltonian. It is natural to expect that the first order corrections can be computed from the normal form obtained in this way. However, due to the quasi-linear nature of the equations, it is not trivial to provide rigorous bounds on the changes of coordinates that allow to construct a reducible (up to certain order) normal form around the torus. Thus we follow a different strategy. \\
To explain it we first compare different Birkhoff normal form approaches
 whose aim is to normalize different type of monomials in the Hamiltonian.
 
 \vspace{0.2em}
\noindent
$\bullet$ \emph{Full} BNF:   \emph{all}
the monomials of degree three and four are normalized;

 \vspace{0.2em}
\noindent
$\bullet$ \emph{Weak + Linear} BNF: it is divided into two steps:
(i) normalization of  all the cubic and quartic terms with at most
one wave number outside the set of tangential sites $S$ (weak); 
(ii) normalization of  all the cubic and quartic terms with exactly 
two wave numbers outside the set  $S$ (linear).
 We remark that both procedures are not convergent.
 We follow a strategy close to the weak+linear BNF and 
 we prove that the corrections found with our approach coincide with the ones we expect.
 We do not explicitly compute such corrections but we 
 prove a uniqueness result on the normal form which implies an a posteriori identification. In non resonant cases this kind of results are provided by classical arguments, see for instance \cite{KdVeKAM}.
For the resonant equation \eqref{eq:113} the main difficulty is due to the presence of non trivial resonances at order four, the Benjamin-Feir resonances, that a priori do not guarantee the integrability nor the uniqueness 
of the normal form around the torus. 
To overcome this problem we would like to obtain the same cancellations on the resonant coefficients of the Hamiltonian obtained in \cite{Zak2}. We do that by exploiting the approximate constants of motion of the system (in the spirit of \cite{FGPa}). 
\\
The main idea is the following: given two Hamiltonians $H=H^{(2)}+O(u^3)$ and $K=K^{(2)}+O(u^3)$ that commute (up to some order) there exists a change of coordinates $\Psi$ that puts simultaneously the Hamiltonians in normal form (up to the same order) namely 
 \[
 H\circ \Psi=H^{(2)}+Z+R, \qquad K\circ\Psi=K^{(2)}+W+Q,
 \]
where the normalized terms $Z, W$ are such that $\{ H^{(2)},  Z\}=0$, $\{ K^{(2)}, W\}=0$ and $R$, $Q$ have higher order of homogeneity. Moreover 
\[
\{ Z, K^{(2)}\}=\{ H^{(2)}, W\}=0\,.
\] Hence the normalized Hamiltonian $Z$ (as well as $W$) is Fourier supported on the common resonances of the adjoint actions of $H^{(2)}$ and $K^{(2)}$. This implies that the coefficients of $Z$ related to non common resonances are automatically zero.\\
We exploit the formal integrability at order four to construct an approximate constant of motion $K$ for the Hamiltonian \eqref{Hamiltonian} with the following property: the adjoint action of $K^{(2)}$ possesses only trivial resonances. This allows to recover \emph{a posteriori} the cancellations of the coefficients of the non-integrable resonant monomials of the normal form and provide the uniqueness argument.

\smallskip
\noindent
 \emph{Invertibility 
of the linearized operator.} 
In a suitable set of coordinates the linear dynamics 
of the tangential (to the torus) variables is decoupled by 
the dynamics of the normal ones (we follow the Berti-Bolle method \cite{BertiBolle}). 
The main difficulty is to invert the linearized operator in the normal directions $\mathcal{L}_{\omega}=\mathcal{L}_{\omega}(U_n)$ at each step of the Nash-Moser scheme. We have that $\mathcal{L}_{\omega}$ is given, up to finite rank operators, by
\begin{equation*}
\omega\cdot\pa_{\vphi}+
\left(
\begin{matrix}
\pa_{x}V+G(\eta)B & -G(\eta) \\
(1+BV_x)+BG(\eta)B & V\pa_{x}-BG(\eta)
\end{matrix}
\right)\,
\end{equation*}
where 
\begin{align*} 
& V =  V (\eta, \psi) :=  (\pa_x \Phi) (x, \eta(x)) = \psi_x - \eta_x B \, , 
\\
& B =  B(\eta, \psi) := (\pa_y \Phi) (x, \eta(x)) =  \frac{G(\eta) \psi + \eta_x \psi_x}{ 1 + \eta_x^2} \, .
\end{align*}
We observe that $\mathcal{L}_{\omega}$ is a pseudo differential operator of order one.
Our aim is to obtain the invertibility of $\mathcal{L}_{\omega}$ and to provide suitable \emph{tame} estimates on the inverse. We do that by means of a reducibility argument, namely we find a quasi-periodically time dependent change of coordinates $\mathcal{T}=\mathcal{T}(\omega t)$ that conjugates $\mathcal{L}_{\omega}$ to a diagonal (in the Fourier basis) operator.
This strategy consists into two main steps:
\begin{itemize}
\item[(a)]\label{stepREGU} a pseudo differential reduction 
in decreasing order of $\mathcal{L}_{\omega}$
which conjugates the linearized operator to 
a pseudo differential  one
with constant coefficients
up to a bounded remainder;

\item[(b)] a reduction of bounded operators and a KAM
 scheme which complete the diagonalization;
 \end{itemize}

 The main new issues are the following:
 \begin{itemize}
 \item[(I)] by the complete resonance the diophantine constants $\gamma$ in 
 \eqref{zerointro} and $\eta$ in  \eqref{secondintro}
 are small with the size of the amplitudes (see for instance 
 \eqref{eps2}). This implies that many terms in $\mathcal{L}_{\omega}$ are not perturbative for the KAM scheme (b), which requires a smallness condition like $\varepsilon^a \gamma^{-1}\ll 1$ for some $a>0$.
 \item[(II)] We consider the linearization on a quasi-periodic traveling embedding $U(\varphi, x)$ without any assumption on the parity of $\varphi$ and $x$. Usually such conditions provide some algebraic cancellations in performing steps $(a)$-$(b)$ and reduce the multiplicity of the eigenvalues simplifying the proof.
 \end{itemize}

We deal with $(\mathrm{I})$ by splitting each step of the reducibility into $2$ parts. In the first part we treat the non perturbative terms of $\mathcal{L}_{\omega}$ by using algebraic arguments to deal with resonances between the linear frequencies of oscillation. We refer to these steps as \emph{preliminary steps} and \emph{linear Birkhoff normal form steps} when they are performed in the procedures described in items $(a)$ and $(b)$ respectively. After the procedures $(a)$ and $(b)$ we obtain a map $\mathcal{T}=\mathcal{T}(\omega t)$ such that 
\[
\mathcal{T} \mathcal{L}_{\omega} \mathcal{T}^{-1}=\omega\cdot \partial_{\varphi}+D, \qquad D=\mathrm{diag}_{j\in\mathbb{Z}} (d_j)
\]
where 
\[
d_j:=\mathfrak{m}_1\,j
+(1+\mathfrak{m}_{\frac{1}{2}})\,\sqrt{|j|}+
\mathfrak{m}_0\,\sign(j)+r_{j}
\]
with $\mathfrak{m}_{1/2}$, $\mathfrak{m}_{0}$, $r_{j}$ depending on $\omega$ and satisfying
\begin{equation}\label{m1}
\mathfrak{m}_1=\varepsilon^2 m_1+O(\varepsilon^4), \quad m_1:=\frac{1}{\pi}\sum_{n\in S} n|n| \zeta_{n}\,,
\end{equation}
\[
(\sup_{j}|j|^{-1/2}|r_{j}|+|\mathfrak{m}_{1/2}|+|\mathfrak{m}_0|)\gamma^{-1}\ll 1\,.
\]
We point out that the first order corrections $O(\varepsilon^2)$ of the $d_j$'s are fundamental to impose the non-resonance conditions since the linear equation \eqref{eq:113CWW} is completely resonant.
Usually such corrections comes from the normalization of the Hamiltonian terms depending on $2$ normal variables. We remark that the weak BNF leaves untouched these terms. Then we use an identification argument to prove that this normalization is obtained through the preliminary steps and the linear Birkhoff normal form aforementioned.  
 
 Concerning $(\mathrm{II})$, we have to deeply exploit the Hamiltonian structure
 and the conservation of momentum. More precisely the map $\mathcal{T}$ constructed in steps $(a)$, $(b)$ is obtained as a composition of several \emph{symplectic} maps that preserve the subspace of traveling embeddings $U(\varphi, x)=\widetilde{U}(\varphi-\mathtt{v} x)$. This allows to obtain algebraic cancellations that provide the diagonalization in decreasing order of the operator $\mathcal{L}_{\omega}$.\\
Now we describe an important step of the pseudo differential reduction of $\mathcal{L}_{\omega}$, namely the reduction of the highest order term where we have to deal with the \emph{singular perturbation} problem. For simplicity we consider the following model example
  \begin{equation}\label{figaro10}
\omega\cdot\pa_{\vphi}+V(\vphi,x)\pa_{x}\,+(1+a(\varphi, x)) |D|^{1/2}.
\end{equation}
The operator $\mathcal{L}_{\omega}$ assumes a similar form after a block-diagonalization and passing to "good unknown" variables. Our aim is to find a change of coordinates that conjugates \eqref{figaro10} to
\begin{equation}\label{figaro11}
\omega\cdot\pa_{\varphi}+\mathfrak{m}_1\,\pa_{x}\,+(1+a_+(\varphi, x)) |D|^{1/2}+\dots
\end{equation}
where $\dots$ denotes lower order pseudo differential operators and $\mathfrak{m}_1$ is in \eqref{m1}.
 Actually this problem is equivalent to straighten the \emph{degenerate} 
 vector field
 $
 \omega\cdot\frac{\pa}{\pa_{\vphi}}+V(\vphi,x)\frac{\pa}{\pa{x}}\,
 $
 on the $\nu+1$-dimensional torus. A priori, due to the degeneracy $V=O(\varepsilon)$, we cannot eliminate the $x$-dependence of the vector field by perturbative arguments. 
 By using the conservation of momentum we can reduce one degree of freedom and study the following \emph{non-degenerate } vector field 
   \[
  (\omega-{\widetilde{V}}(\Theta)\mathtt{v})\cdot\frac{\pa}{\pa{\Theta}}\,,
  \qquad \Theta\in \mathbb{T}^{\nu}\,,
  \]
  where $\widetilde{V}(\vphi-\mathtt{v}x)=V(\vphi,x)$. Then we can apply a result of straightening of weakly perturbed constant vector fields on tori given in \cite{FGMP}.  We remark that $\mathfrak{m}_1\approx (2\pi)^{-\nu}\int_{\mathbb{T}^{\nu}} \widetilde{V}(\Theta)d\Theta $.\\
 When we conjugate the operator \eqref{figaro10} into \eqref{figaro11} we need to ensure that the pseudo differential structure is preserved and we need to provide suitable estimates for the lower order terms. We do that by using a quantitative Egorov theorem developed in \cite{FGP}.
  
  \smallskip
\noindent
$(iv)$ 
  \emph{Non-resonance conditions.}
  Since the dispersion law is \emph{sub-linear}, 
 at each step of 
  items $(a)$ and $(b)$ small divisors problems arise. We discuss 
  the non-resonance conditions that we shall require on the eigenvalues.
In order to deal with the operator in \eqref{figaro10}
we shall impose the following
``zero order Melnikov conditions'':
\begin{equation}\label{zerointro}
|\big(\omega-\mathfrak{m}_1\mathtt{v}\big)\cdot\ell|
\geq \gamma \langle \ell\rangle^{-\tau}\,, \qquad \forall\, \ell\in \mathbb{Z}^{\nu}\setminus\{0\}.
\end{equation}

In the KAM reducibility scheme in item $(b)$
we shall  impose suitable lower bounds on the function
\begin{equation*}
\psi_{\ell, j,k}=\omega\cdot\ell +d_j-d_k
\end{equation*}
with $\ell\in \mathbb{Z}^{\nu}$, $j,k\in \mathbb{Z}\setminus\{0\}$.
These conditions are called ``second order Melnikov conditions'' and are of the form
\begin{equation}\label{secondintro}
|\psi_{\ell,j,k}|\geq \frac{\eta}{\langle \ell\rangle^{\tau}}\,,\quad
\ell\in\mathbb{Z}^{\nu}\,,\;\; j,k\in\mathbb{Z}\setminus\{0\}\,,\quad (\ell,j,k)\neq (0,j,j)\,,
\end{equation}
where $\eta, \tau>0$ are some constants to be fixed.
We prove that
\eqref{secondintro} hold with
\begin{equation}\label{secondintrobis}
\eta:=\gamma^{3}\ll \gamma\,,
\end{equation}
where $\gamma$ is the diophantine constant in \eqref{zerointro}.
Actually we show that  the complementary of the set of  frequencies $\omega$
such that \eqref{secondintro} hold has measure
going to zero as $\e\to 0$.
The crucial problem is the summability in the indexes $\ell,j,k$.
Since in this case the linear 
frequencies grow at infinity in sublinear way, 
the difference $\sqrt{|j|}-\sqrt{|k|}$ accumulates everywhere in $\mathbb{R}$.
This means that for any fixed $\ell$ 
there are \emph{infinitely many} indexes $j,k$ to be taken into account. 
The key idea that we use is the following. By the conservation of momentum
we 
 need to impose the conditions
\eqref{secondintro} only for $(\ell,j,k)$ satisfying
\begin{equation*}
\mathtt{v}\cdot\ell+j-k=0\,,
\end{equation*}
where $\mathtt{v}$ is the velocity vector given in \eqref{velocityvec}.
This allows us to show that, if $|j|,|k|$
are much larger than $|\ell|$, then
the conditions \eqref{secondintro} are implied by the
\eqref{zerointro}
and then we are left 
 to control the small divisors only for \emph{finitely many} 
indexes $j,k$.

An alternative approach would be to prove 
Melnikov conditions 
by setting
\begin{equation}\label{muloss}
\eta:=\frac{\gamma}{\langle j\rangle^{\mathtt{d}}\langle k\rangle^{\mathtt{d}}}\,,
\quad \mathtt{d}> 1\,.
\end{equation}
This choice would allow to prove that ``many'' frequencies $\omega$ satisfy
the \eqref{secondintro},
but with the disadvantage of having small divisors that create a loss 
of space derivatives.

The class of non perturbative terms in the case $\gamma^3=\eta$ 
is larger respect to the choice $\gamma=\eta$, 
however such terms are compactly time-Fourier 
supported ($|\ell|\leq C$ for some constant $C>0$).
In \cite{FGww}, in order to deal with such non perturbative terms, we 
adopt the following strategy:
if at least one between $|j|, |k|$
is large enough
then there are no small divisors; 
otherwise we choose to impose 
the second order resonance Melnikov conditions 
with $\eta$ as in \eqref{muloss}.
Actually this does not create any loss of 
derivatives since $|j|,|k|$ are taken into a ball with finite radius.
We deal with the perturbative terms by imposing conditions \eqref{secondintro}
with $\eta$ as in \eqref{secondintrobis}.


\bigskip
\begin{thebibliography}{12}

\bibitem{BBHM}
P. Baldi, M. Berti, E. Haus, R.  Montalto, 
\emph{Time quasi-periodic gravity water waves in finite depth},
Invent. math. 214, 739--911, (2018). 


\bibitem{BertiBolle} M. Berti, P. Bolle, 
\textit{A Nash-Moser approach to KAM theory}, 
Fields Institute Communications, volume 
75 ``Hamiltonian PDEs and Applications", 
255-284, (2014).



\bibitem{BM1}
M. Berti, R.  Montalto, {\it KAM for gravity capillary water waves}, 
{Memoires of AMS}, Memo 891, vol. 263 (2020).

\bibitem{CN1}
W. Craig, D. Nicholls, \emph{Travelling two and three 
dimensional capillary gravity water waves}, SIAM J. Math. Anal.,
32(2):323-359, (2000).

\bibitem{CrSu}
W. Craig, C. Sulem,
\newblock {\it Numerical simulation of gravity waves},
\newblock { J. Comput. Phys.}, 108(1):73--83, (1993).

\bibitem{CW}
W. Craig, P. Worfolk, 
\newblock {\it An integrable normal form for water waves in infinite depth},
\newblock { Phys. D}, 84(3-4):513-531, (1995). 

\bibitem{Zak2}
A. I. Dyachenko, V. E. Zakharov,   {\it Is free-surface hydrodynamics an integrable system?}
Physics Letters A, 190(2):144--148, (1994).

\bibitem{FGww} R. Feola, F. Giuliani, 
\emph{Quasi-periodic Traveling Waves on an Infinitely Deep Perfect Fluid Under Gravity},
preprint, arXiv:2005.08280 (2020).

\bibitem{FGMP} R. Feola, F.  Giuliani, R.  Montalto, M. Procesi, 
\textit{Reducibility of first order linear operators on tori via Moser's theorem}, 
Journal of Functional Analysis
276(3):932-970 (2019). 
DOI: https://doi.org/10.1016/j.jfa.2018.10.009.

\bibitem{FGPa} R. Feola, F. Giuliani, S. Pasquali, 
\textit{On the integrability of Degasperis-Procesi equation: control of the 
Sobolev norms and Birkhoff resonances}, 	
Journal of Differential Equations 266(6):3390--3437, (2019).
DOI: https://doi.org/10.1016/j.jde.2018.09.003. 



\bibitem{FGP}
R. Feola, F.  Giuliani, M. Procesi, \emph{Reducibility for a class 
of weakly dispersive linear operators arising from 
the Degasperis Procesi equation},
Dynamics of Partial Differential Equations, 16(1): 25-94 (2019).
DOI: http://dx.doi.org/10.4310/DPDE.2019.v16.n1.a2.

\bibitem{FGP1}
R. Feola, F. Giuliani, M. Procesi, \emph{Reducibile KAM tori for the 
Degasperis Procesi equation}, Comm. Math. Phys. (2020), 
 DOI: 10.1007/s00220-020-03788-z.


\bibitem{IossPlo111} G. Iooss, P.  Plotnikov, 
\emph{Multimodal standing gravity waves: 
a completely resonant system}, 
J. math. fluid mech. 7, S110--S126, (2005). 

\bibitem{IoossPlo1} G. Iooss, P. Plotnikov, \emph{Small divisor problem 
in the theory of three-dimensional water gravity waves}, 
Mem. Amer. Math.
Soc., 200(940), (2009).

\bibitem{IoossPlo2} G. Iooss, P. Plotnikov, \emph{Asymmetrical tridimensional 
traveling gravity waves}, Arch. Rat. Mech. Anal., 200(3):789--880,
(2011).



\bibitem{IoossPloTol1} G. Iooss, P. Plotnikov, J. Toland, 
\emph{Standing waves on an infinitely deep perfect fluid under gravity}, 
Arch. Ration. Mech.
Anal., 177(3):367--478, 2005.

\bibitem{KdVeKAM} T. Kappeler, J. P\"oschel, \textit{KAM and KdV}, Springer (2003).

\bibitem{PloTol1} P. Plotnikov, J. Toland, \emph{Nash-Moser theory for 
standing water waves}. Arch. Ration. Mech. Anal., 159(1):1-83, 2001.


\bibitem{Zak1}  V.E. Zakharov, {\it Stability of periodic waves of 
finite amplitude on the surface of a deep fluid}, 
J Appl Mech Tech Phys 9:190--194, (1968). 




\end{thebibliography}
\end{document}